\documentclass[11pt,letterpaper]{article}
\usepackage{amsmath,latexsym,amssymb,ifthen}

\renewcommand{\emph}[1]{\textit{#1}}
\usepackage{enumerate}
\usepackage{color}\usepackage{graphicx}

\definecolor{brown}{cmyk}{0, 0.72, 1, 0.45}
\definecolor{grey}{gray}{0.5}

\newcommand{\ignore}[1]{}

\newboolean{draft}
\setboolean{draft}{true}
\def\cA{{\cal A}}
\def\cB{{\cal B}}

\def\cT{\mathcal{T}}

\newcommand{\set}[1]{\left\{#1\right\}}

\newcommand{\proofstart}{{\noindent \bf Proof\hspace{2em}}}
\newcommand{\proofend}{\hspace*{\fill}\mbox{$\Box$}}

\def\ii_(#1,#2){i_{#1}^{#2}}

\def\LL{\Lambda}

\def\a{\alpha}
\def\b{\beta}
\def\d{\delta}

\def\e{\varepsilon}
\def\f{\phi}
\def\F{\Phi}
\def\g{\gamma}
\def\G{\Gamma}

\def\Th{\Theta}
\def\l{\lambda}

\def\n{\nu}
\def\p{\pi}

\def\r{\rho}

\def\s{\sigma}
\def\t{\tau}
\def\om{\omega}

\def\cW{{\cal W}}

\newcommand{\rdup}[1]{\left\lceil #1 \right\rceil}
\newcommand{\rdown}[1]{\mbox{$\left\lfloor #1 \right\rfloor$}}

\newcommand{\ooi}{(1+o(1))}

\def\cE{\mathcal{E}}

\newcommand{\brac}[1]{\left( #1 \right)}

\newcommand{\expect}{\operatorname{\bf E}}
\def\E{\expect}
\renewcommand{\Pr}{\operatorname{\bf Pr}}
\newcommand\bfrac[2]{\left(\frac{#1}{#2}\right)}

\newtheorem{theorem}{Theorem}
\newtheorem{lemma}{Lemma}
\newtheorem{corollary}[lemma]{Corollary}

\newtheorem{remthm}[lemma]{Remark}

\newenvironment{remark}{\begin{remthm}\rm }{\end{remthm}}%
\newcounter{thmtemp}

\newcommand{\nospace}[1]{}

\def\path{\operatorname{PATH}}

\newcommand{\beq}[1]{\begin{equation}\label{#1}}
\def\eeq{\end{equation}}

\newcommand{\card}[1]{\left|#1\right|}

\def\bt{{\bf t}}

\def\cR{{\cal R}}

\def\integer{\mathbb{Z}}
\def\cX{{\cal X}}

\newcommand{\bP}[1]{{\bf P#1}}

\newcounter{rot}%\addtocounter{rot}{1}, \therot

\begin{document}

\title{A note on the vacant set of  random walks on the hypercube and other regular graphs of high degree}
\author{Colin Cooper\thanks{Department of  Informatics,
King's College, University of London, London WC2R 2LS, UK.
Research supported  by EPSRC grant EP/J006300/1.}
\and Alan
Frieze\thanks{Department of Mathematical Sciences, Carnegie Mellon
University, Pittsburgh PA 15213, USA.
Supported in part by NSF grant CCF1013110.
}}
\date{\today}
\maketitle
\begin{abstract}
We consider a random walk on a $d$-regular graph $G$ where $d\to\infty$
and $G$ satisfies certain conditions. Our prime example is the
$d$-dimensional hypercube, which has $n=2^d$ vertices. We
explore the likely component structure of the
vacant set, i.e. the set of unvisited vertices. Let $\LL(t)$ be the subgraph induced by the vacant set
of the walk at step $t$.
We show that if certain conditions are satisfied then the graph $\LL(t)$ undergoes
 a phase transition at around $t^*=n\log_ed$.
Our results are that if $t\leq(1-\e)t^*$ then w.h.p. as the number
vertices $n\to\infty$, the size $L_1(t)$ of the largest component
satisfies $L_1\gg\log n$ whereas if $t\geq(1+\e)t^*$ then
$L_1(t)=o(\log n)$.
\end{abstract}
\section{Introduction}
The problem we consider can be described as follows.
We have a finite graph $G=(V,E)$, and a simple random walk $\cW=\cW_u$ on $G$, starting at $u\in V$. In this walk, if $\cW(t)$ denotes the position of the walk after $t$ steps, then $\cW(0)=u$ and if $\cW(t)=v$ then $\cW(t+1)$ is equally likely to be any neighbour of $v$. We consider the likely component structure of the subgraph $\LL(t)$ induced by the unvisited vertices of $G$ at step $t$ of the walk.

Initially all vertices $V$ of $G$ are unvisited or {\em vacant}.
We regard unvisited vertices as colored {\em red}.
When $\cW_u$ visits a vertex, the vertex is re-colored {\em blue}.
Let $\cW_u(t)$ denote the position of $\cW_u$ at step $t$.
Let $\cB_u(t)=\set{\cW_u(0),\cW_u(1),\ldots,\cW_u(t)}$ be the set of
blue vertices at the end of step $t$, and $\cR_u(t)=V\setminus \cB_u(t)$.
Let $\LL_u(t)=G[\cR_u(t)]$ be the subgraph of $G$ induced by $\cR_u(t)$.
Initially $\LL_u(0)$ is connected, unless $u$ is a cut-vertex. As
the walk continues, $\LL_u(t)$ will shrink to the empty graph
once every vertex has been visited. We wish to determine, as far as
possible, the likely evolution of the component
structure as $t$ increases.

For several graph models, it has been shown that the component structure of $\LL(t)=\LL_u(t)$
undergoes a phase transition of some sort. In this paper we add results for some new classes of graphs.
What we expect to happen is that there is a time $t^*$, such that if $t\geq (1+\e)t^*$ then w.h.p. all components of $\LL(t)$ are ``small'' and if $t\leq (1-\e)t^*$ then w.h.p. $\LL(t)$ contains some ``large'' components. Here $\e$ is some arbitrarily small positive constant and the meanings of small, large will be made clear.

\subsection{Previous work}

We begin with the paper by \v{C}ern\'y, Teixeira and Windisch \cite{CTW}. They consider a sequence of $n$-vertex graphs $G_n$ with the following properties:
\begin{enumerate}
\item[A1] $G_n$ is $d$-regular, $3\leq d=O(1)$.
\item[A2] For any $v\in V(G_n)$, there is at most one cycle within distance $\a\log_{d-1}n$ of $v$ for some $\a\in (0,1)$.
\item[A3] The second eigenvalue $\l_2$ of the random walk transition matrix satisfies $\l_2\leq 1-\b$ for some constant $\b\in (0,1)$.
\end{enumerate}
Let
\beq{t*}
t^*=\frac{d(d-1)\log(d-1)}{(d-2)^2}n.
\eeq
In which case, it is shown in \cite{CTW} that for $t\leq (1-\e)t^*$ there is w.h.p. a unique giant component in $\LL(t)$ of size $\Omega(n)$ and other components are all of size $o(n)$. Furthermore, if $t\geq (1+\e)t^*$ then all components of $\LL(t)$ are of size $O(\log n)$.

The most natural class of graphs satisfying A1,A2,A3 are random $d$-regular graphs, $3\leq d=O(1)$. For this class of graphs Cooper and Frieze \cite{CF1} tightened the above results in the following ways.  (i) they established the asymptotic size of the giant component for $t\leq (1-\e)t^*$, and proved all other components have size $O(\log n)$; (ii) they proved almost  all small components are trees, and gave a detailed census of the number of trees of sizes $O(\log n)$. Subsequent to this work, \v{C}erny and Teixeira \cite{CT} built on the methodology of \cite{CF1} and analysed the component structure at time $t^*$ itself.
More recently, for random $d$-regular graphs, $3\leq d=O(1)$, Cooper and Frieze \cite{CFVacNet} determined   the phase transition for a related structure, the {\em vacant net}, which by analogy with vacant set, they define as the subgraph
induced by the unvisited edges of the graph $G$. Initially all edges are unvisited. The random walk {\em visits an edge} by making a transition using the edge.

In the  paper \cite{CF1}, Cooper and Frieze  also considered the class of Erd\H{o}s-Re\'nyi random graphs  $G_{n,p}$
with edge probabilities $p$ above the connectivity threshold $p={\log n}/{n}$.
For $G_{n,p}$ where $p={c\log n}/{n}$, $(c-1)\log n\to \infty$, they established that $\LL(t)$ undergoes a phase transition  around $t^*=n\log\log n$. For these graphs, at $t_{-\e}=(1-\e)t^*$ the size $L_1$ of the largest component cannot be $\Omega(n)$ since the vacant set has size $|\cR(t_{\e})|=o(n)$ w.h.p.  On the other hand it was shown that $L_1=\Omega(|\cR(t_{\e})|)$ w.h.p. More recently, Wassmer \cite{W} found the phase transition in $\LL(t)$ when the underlying graph is the giant component of $G_{n,p}$, $p=c/n$, $c>1$.

There has also been a considerable amount of research on the $d$-dimensional grid $\integer^d$ and the $d$-dimensional torus $(\integer/n\integer)^d$. Here the results are less precise. Benjamini and Sznitman \cite{BS} and  Windisch \cite{DW}
investigated the structure of the vacant set of a random walk
on a $d$-dimensional torus. The main focus of this work is to apply the method of
random interlacements. For toroidal grids of dimension $d \ge 5$, it is shown that
there is a value $t^+(d)$,  linear in $n$, above which the vacant set is sub-critical, and a value of
$t^-(d)$ below which the graph is super-critical. It is believed that there is
 a phase transition for $d \ge 3$.
A recent monograph by \v{C}erny and Teixeira \cite{CTm}
summarizes the random interlacement methodology. The monograph also gives details for the
vacant set of random $r$-regular graphs.

\subsection{New results}
In this note we consider certain types of $d$-regular graphs with $n$ vertices, where $d\to\infty$ with $n$.
Our main example of interest is the hypercube $Q_d$ which has $n=2^d$ vertices. The vertex set of the hypercube is sequences $\set{0,1}^d$ where two sequences are defined as adjacent iff they differ in exactly one coordinate. In order to be slightly more general, we identify those properties of the hypercube that underpin our results.

Given certain properties (listed below), we can show that w.h.p. the graph $\LL(t)$ exhibits a change in component structure at around the time $t^*=n\log d$ which is asymptotically equal to the expression in \eqref{t*}. We show that if $t\leq t_{-\e}=(1-\e)t^*$ then w.h.p. there are components in $\LL(t)$ of size much larger than $\log n$, whereas if $t\geq t_\e=(1+\e)t^*$ then all components of $\LL(t)$ are of size $o(\log n)$.

We use the notation $\Pr(\cW_x(t)=y)$ and  $P^{t}_x(y)$ for the probability that a random walk starting from vertex $x$ is at vertex $y$ at step $t$. If $t$ is sufficiently large, so that the walk is very close to stationarity and the starting point $x$ is arbitrary, we may also use
the simplified notation $\Pr(\cW(t)=y)$.
 Let $\pi_v=d(v)/2m$ to denote the stationary probability of vertex $v$, where $m=|E|$ is the number of edges of the graph $G$ and $d(v)$ is the degree of $v$. For regular graphs, $\pi_v=1/n$.
The rate of convergence of the walk is given by
\begin{equation}\label{mix}
|P^{t}_x(y)-\pi_y| \leq (\p_y/\p_x)^{1/2}\l^t,
\end{equation}
where $\l=\max(\l_2,|\l_n|)$ is the second largest eigenvalue of the transition matrix
in absolute value.
See for example,  Lovasz \cite{Lo} Theorem 5.1.

The  hypercube $Q_d$ is bipartite, and a simple random walk  does not have a stationary distribution on  bipartite graphs.
To overcome this, we  replace the simple random walk by a {\em lazy} walk,
in which at each step there is a 1/2 probability of staying put.
Let $N_G(v)$ denote the neighbours of $v$ in $G$, and $d_G(v)=|N_G(v)|$. The lazy walk $\cW$ has transition probabilities $P^{t}_v(w)$ given by
\[
P^{t}_v(w)=\begin{cases}\frac12&w=v\\ \frac1{2 d_G(v)}&w\in N_G(v)\\0&\text{Otherwise}\end{cases}.
\]
We can obtain the underlying simple random walk, which we refer to as the {\em speedy} walk,
 by ignoring the steps when the particle does not move. For large $t$, asymptotically half of the steps in the lazy walk will not result in a change of vertex. Therefore  w.h.p. properties of the speedy walk after approximately $t$ steps can be obtained from properties of the lazy walk after approximately $2t$ steps.

The effect of making the walk lazy is to shift the eigenvalues of the simple random walk upwards so that, for the lazy walk $\l=\l_2$.
We define a mixing time $T$ for the lazy walk by
\beq{4}
T=\min_{t\geq 1}\set{t:\;\card{P^{t}_x(y)-\frac1n}\leq \frac1{n^3}}.
\eeq
 For the lazy walk, the spectral gap is $1-\l$, so using this in \eqref{mix}, property \bP1
 implies that we can take $T=O(d^{\r_1}\log n)$ in \eqref{4}.
  %(see \cite{LPW}, Chapter 12).

\subsubsection*{The graph properties we assume for our analysis}

Let $G=(V,E)$ be a graph with vertex set $V$ and edge set $E$.
For  $S \subset V$, define $N_G(S)=\set{w\in V\setminus S:\;\exists v\in S\ s.t.\  \set{v,w}\in E}$.

We assume that the graph $G=(V,E)$ is $d$-regular, connected, and has the properties \bP1--\bP4 listed below.
The bounds in  properties \bP2--\bP4 are parameterised by the $\e$ used to
define $t_{\pm\e}$ for the vacant set. We will point out later where we use these bounds,
so that the reader can see their relevance.

\begin{enumerate}[{\bf P1}]
\item The spectral gap for the lazy walk is $\Omega(1/d^{\r_1})$ for some
  constant $0<\r_1 \le 3$.
 This  implies that we can take $T=O(d^{\r_1}\log n)$ in \eqref{4},
  (see \cite{LPW}, Chapter 12).

  \item $(\log\log n)^{2/\e}\ll d =O\bfrac{n}{\log n}^{1/5}$.
\item For $u,v\in V$, the graph distance $dist_G(u,v)$ is the length of the shortest path from $u$ to $v$ in $G$.
Let $\n(u,v)$ be the number of neighbours $w$ of $v$ for which $dist_G(w,u)\leq dist_G(u,v)$.  Let $\r_2=O(1)$.
Then for all $u,v$ such that $dist_G(u,v)\leq d^\e$, we have $\n(u,v)\leq \r_2\;dist_G(u,v)$.
%\item There exists $\a=\a(d)>0$ such that $\a d^\e\gg \log d$, and for any  $S\subseteq V,|S|\leq \frac{\a n}{d}$ the inequality  $|N_G(S)|\geq \a d|S|$ holds. We do not need $\a$ to be a constant here. If $\a\to 0$, then \bP4 can still hold provided $d$ grows fast enough.
\item  For $S\subseteq V$, let $e(S)$ denote the number of edges induced by $S$. If  $|S|=o\brac{\log n}$, then $e(S)=o(d|S|)$.
\end{enumerate}
Properties \bP1--\bP4 are various measures of expansion.
Random regular graphs with degree $d$ satisfying \bP2 satisfy the other properties w.h.p. The hypercube $Q_d$ satisfies these properties to a degree.
%\subsubsection*{The theorems of this paper}
Our results for the structure of the vacant set $\LL(t)$ based on these properties are as follows.

\begin{theorem}\label{th1}
Let $\e=\e(n)$ be a function such that  $\e\gg 1/\log d$.
Let $t^*=n\log d$ and let $t_{\pm\e}=(1\pm\e)t^*$. Let $L_1(t)$ denote the size of the largest component in $\LL(t)$.
At step $t$ of the speedy walk, the following results for $L_1(t)$  hold.
\begin{enumerate}[{\bf(a)}]
\item  If $G$ satisfies  \bP1, \bP2, \bP3, \bP4, and  $t\leq t_{-\e}$ then w.h.p. $L_1(t)\geq e^{\Omega(d^{\e/2})}$.\\
Note that $d^{\e/2}$ can be replaced by $d^{\g\e}$ for any constant $0<\g<1$.
\item If $G$ satisfies \bP1, \bP2, \bP3, and
 $t\geq t_{+\e}$ then w.h.p. $L_1(t)=o(\log n)$.
%\item If $G$ satisfies \bP1, \bP2, \bP3, \bP4, and   $t\leq t_{-\e}$ then w.h.p. $L_1(t_{-\e})\geq \frac{\a n}{d}$.
\end{enumerate}
\end{theorem}

We next prove that the hypercube $Q_d$ satisfies Theorem \ref{th1}(a),(b).
To show this, we check that $Q_d$ satisfies properties \bP1--\bP4.
\bP1 is satisfied with $\r_1=1$, as the spectral gap for the lazy walk is
$\frac{2}{d}$ (see \cite{LPW} page 162).
\bP2 is clearly
satisfied, as $d=\log_2n$.
For \bP3, without loss of generality, let $v=(0,0,\ldots,0)$ and let $u=(1,1,\ldots,1,0,\ldots,0)$
($k$ 1's) be vertices of $Q_d$. There are exactly $\n(u,v)=k$ neighbours $w$ of $v$
which satisfy $dist_G(u,w) \le dist_G(u,v)$, so we can take $\r_2=1$.
For \bP4 we can use the edge isoperimetric inequality of Hart \cite{Hart} which
states that the number of edges between $S$ and $V-S$ is at least $s(d-\log_2 s)$, where $|S|=s$.
This implies that $S$ induces at most $(s/2)\log_2 s$ edges.
If $s=o(d)$ then $e(S)\leq (s/2)\log_2s=o(ds)$.
%We do not know if the condition \bP4 holds for $Q_d$. The value of $\a$  obtained from the vertex isoperimetric inequality for $Q_d$  is too small.

\section{The main tools for our proofs}

Given a graph $G$ and random walk $\cW$, let $T$ be the mixing time given in \eqref{4}.
 For a vertex $v$, let $R_v=R_v(G)$ denote the expected number of visits to $v$ by the speedy walk $\cW_v$ within $T$ steps. Thus
 \beq{RvDef}
 R_v=\sum_{k=0}^{T} P_v^k(v).
 \eeq
Note that, as $P_v^0(v)=1$,  $R_v\geq 1$.

Our main tool is a lemma (Lemma \ref{First}) that we have found very useful in analysing the cover time of various classes of random graphs.
A  general form of Lemma \ref{First}  as proved in \cite{CFgiant}  requires a certain technical condition to be satisfied. It was shown in \cite{CFR} that provided $R_v=O(1)$ for all $v\in V$, this condition is always true. We prove in Lemma \ref{lemRv} that if \bP2 and \bP3  hold, then $R_v= 2+O(1/d)$ for all $v \in V$. The probabilities given in Lemma \ref{First} and Corollary \ref{geom} are with respect to the random walk.

\begin{lemma}[First visit lemma]
\label{First}
Suppose that $R_v=O(1)$ for $v\in V$ and $T\pi_v=o(1)$ and $T\pi_v=\Omega(n^{-2})$.
Let
$$f_t(u,v)=\Pr(t=\min\set{\t>T:\cW_u(\t)=v})$$
be the probability that the first visit to $v$ after time $T$
occurs at step $t$.

There exists
\begin{equation}\label{pv}
p_v=\frac{\pi_v}{R_v(1+O(T\pi_v))},
\end{equation}
 and constant $K>0$ such that for all $t\geq T$,
\beq{frat}
f_t(u,v)=(1+O(T\pi_v))\frac{p_v}{(1+p_v)^{t+1}}+O(T \pi_v e^{- t/KT}).
\eeq
\end{lemma}
\proofend

\begin{corollary}
\label{geom}
For $t\geq T$ let $\cA_v(t)$ be the event that $\cW_u$ does not visit $v$ at steps $T,T+1,\ldots,t$.
Then, under the assumptions of Lemma \ref{First},
\beq{atv}
\Pr(\cA_v(t))=\frac{(1+O(T\p_v))}{(1+p_v)^{t}} +O(T^2\p_ve^{- t/KT}).
\eeq
\end{corollary}
The result \eqref{atv} follows by adding up \eqref{frat} for $s \ge t$.
\proofend

\begin{remark}\label{remlab}
Provided $p_v=o(1/T)$ and $t \ge L$ where
\beq{Lval}
L = 2KT \log n
\eeq
then, as $p_v=O(\pi_v)$, bounds \eqref{frat} and \eqref{atv} can be written as
$$f_t(u,v)=(1+O(T\pi_v))\; {p_v}\;{(1-p_v)^{t}}$$
and
$$\Pr(\cA_v(t))={(1+O(T\p_v))}\; {(1-p_v)^{t}}$$
respectively.
For the graphs we consider $\pi_v=1/n$.
From \bP1, $T=O(d^{\r_1} \log n)$ and from \bP2, $d=O(n/\log n)^{1/4}$.
Thus for $\r_1 \le 3$, $p_v=o(1/T)$ as required.
\end{remark}

\subsection*{Contraction lemma}
Let $H=(V(H),E(H))$ be given.
Let $S$ be a subset of vertices of $H$. By contracting $S$ to single vertex
$\g(S)$, we form a multi-graph $\G=\G(H,S)=(V',E')$  in which the set $S$ is replaced by $\g$. The edges of $H$, including loops and multiple edges formed by contraction, are retained. Thus if $(v,w) \in E(H)$ and $v,w \not \in S$ then $(v,w) \in E'$,
 whereas if $v \in S$ and $(v,w) \in E(H)$
then $(\g, w) \in E'$. This includes the case $w \in S$ so that $(\g,\g)\in E'$.
It follows that $|E'|=|E(H)|$.
We prove in Lemma \ref{contract} that the probability of a first visit to $S$  in $H$ can be found (up to an additive error of $O(|S|/n^3)$ from
the probability of a first visit to $\g$ in $\G$.

Note  that if $T$ is a mixing time for $\cW$ in $H$,
 %obtained using the {\em  conductance} of $H$,as is the case in this paper,
then $T$ is a mixing time for the walk in $\G$.
We used the second eigenvalue $\l_2(H)=\l$ of the lazy walk
in \eqref{mix} to obtain
the mixing time bound in \eqref{4}.
It is proved in \cite[Ch.\ 3]{AlFi}, Corollary 27,
that if a subset $S$ of vertices is contracted to a single vertex,
then the second eigenvalue of the transition matrix cannot increase.
Thus $\l_2(H) \ge \l_2(\G)$.

\begin{lemma}\label{contract}\cite{CFgiant}
 Let $H=(V(H),E(H))$, let $S \subseteq V(H)$, let $\g(S)$ be the contraction
of $S$, preserving edges, including loops and multiple edges. Let $V'=V-S+\g$, and let $\G(H)=(V',E')$.
Let $\cW_u$ be a random walk in $H$ starting at $u \not \in S$, and let $\cX_u$ be a random walk in
$\G$. Let $T$ be a mixing time satisfying \eqref{4} in both $H$ and $ \G$.
For graphs $G=H, \G$, let $\cA^G_w(t)$ be the event that in graph $G$,
no visit was made to $w$ at any step $T \le s \le t$. Then
\[
\Pr( \wedge_{v \in S}\; \cA_{v}^H(t)) =\Pr(\cA_{\g}^\G(t))+O(|S|/n^3).
\]
\end{lemma}
\proofstart

Note that $|E(H)|=|E'|=m$, say. Let $W_x(j)$ (resp. $X_x(j)$)  be the position of walk $\cW_x$ (resp. $\cX_x(j)$)
 at step $j$.
For graphs $G=H, \G$, let $P_u^s(x;G)$ be the $s$ step transition probability  for the corresponding walk  to go from $u$ to $x$ in $G$.
\begin{eqnarray}
\Pr(\cA_{\g}^\G(t))&=&\sum_{x\ne \g}P^{T}_{u}(x;\G)
\Pr(X_x(s-T)\neq \g,\; T\leq s\leq t; \G)
\label{whatis}\\
&=&\sum_{x\ne \g}\brac{\frac{d(x)}{2m}+O\brac{1/n^{3}}  }
\Pr(X_x(s-T)\neq \g,\; T\leq s\leq t; \G) \label{close}\\
&=& \sum_{x \not \in S}
\brac{P^{T}_{u}(x;H)+O\brac{1/n^{3}}} \Pr(W_x(s-T)\not \in  S,\;T\leq s\leq t;H)\label{extra1}\\
&=& \sum_{x \not \in S}\brac{\Pr(W_u(T)=x)\Pr(W_x(s-T)\not \in  S,\;T\leq s\leq t;H)+O(n^{-3})}\nonumber\\
&=& \Pr(W_u(t)\not \in S,\;T\leq s\leq t;H)+O(|S|/n^{3})\nonumber\\
&=& \Pr(\wedge_{v \in S}\; \cA_{v}^H(t)) +O(|S|/n^3)
.\label{40}
\end{eqnarray}
In \eqref{whatis}, if $\cA_{\g}^\G(t)$ occurs then $X_u(T) \ne \g$. Given $X_u(T)=x$, by the Markov property $X_u(s)$ is equivalent to the walk $X_x(s-T)$.
 After $T$ steps, the walk $X_u$
on $\G$ is close to stationarity. We use \eqref{4} to approximate $P_u^T(x;\G)$ by $\pi_x=d(x)/2m=1/n$ in \eqref{close}.
The second factor in equation \eqref{extra1} follows because there is a natural measure preserving
map $\f$ between walks in $H$ that start at $x\not \in S$ and avoid $S$, and walks
in $\G$ that start at $x \ne \g$ and  avoid $\g$.
\proofend

\begin{remark}
We use Lemma \ref{contract}  throughout the rest of this paper, and often without further mention. Indeed most of the proofs rely on contracting some set $S$ or other to a vertex $\g(S)$. In this case, although a different graph $\G$, and different walk $\cX$ are used to estimate the probability, provided
\[
\frac{|S|}{n^3} = o(\Pr(\cA_{\g}^\G(t)) ),
\]
the probability estimate we obtain for the walk $\cW$ in the base graph $H$ is correct. If necessary, by increasing the mixing time $T$ by a constant factor we can reduce the error term $|S|/n^3$ to $|S|/n^c$ for any $c>0$.
\end{remark}

\subsection*{Visits to sets of vertices}

Given   the walk made a first visit  to set of vertices $S$, we  need  the
probability  this first visit was  to a given $v \in S$.

\begin{lemma}\label{whichvx}
Let $S=\{v_1,...,v_k\}$ be a set of vertices of a graph $G$, such that
 Lemma \ref{First} holds in $G$ for all $v \in S$, and  also for $\g(S)$ in $\G(G)$.
For $t \ge T$, let $ \cE_v=\cE_v(t)$ be the event that the first visit to $v$
after time $T$ occurs at step $t$,
(i.e. $t=\min\set{\t>T:\cW(\t)=v}$),
and let $\cE_S=\cup_{v \in S} \cE_v$.
Suppose $t\ge 2(T+L)$ where $L= 2KT \log n$, where $K>0$ is some suitably large constant then for $v \in S$
\beq{myvertex}
\Pr(\cE_v \, \mid \, \cE_S)= \frac{p_v}{\sum_{w \in S} p_w} (1+O(\xi)),
\eeq
where $\xi=L \pi_S$, and $p_w, w \in S$ are as defined in Lemma \ref{First} for the walk on $G$.
\end{lemma}
\proofstart
It is enough to prove the lemma for $S=\{u,v\}$, i.e. for two vertices,  as vertex $u$
can always be a contraction of a set. Specifically, if $|S|>2$ let $u=\g(S\setminus\{v\})$.

Write $t$ as $t=T+s+T+L$, where $s \ge L$.
Let $\cA_u$ be the event that  $\cW(t)=u$, but that $\cW(\s) \not \in \{u,v\}$ for  $ \s \in[T,s+T-1]$,  and that $\cW(\s) \ne u$ for $ \s \in [s+2T,t-1]$.
Contract $S$ to  $\g(S)$ and apply Corollary \ref{geom}, Remark \ref{remlab}
  and Lemma \ref{contract} to $\g(S)$ in $[T,T+s-1]$. The probability of no visit to $S$ is $(1+O(T\pi_S)) (1-(p_u+p_v))^s$.
Next, apply Lemma \ref{First} (and Remark \ref{remlab}) to $u$ in $[s+2T,t]=[t-L,t]$.
The probability of a first visit to $u$ is $(1+O(T\pi_u))(1-p_u)^{L}p_u$. Thus
\beq{ohno}
\Pr(\cA_u) \le (1+O(T\pi_S)) (1-(p_u+p_v))^s\;(1-p_u)^{L}p_u.
\eeq

Let $\cB_u$ be the event that $\cW(t)=u$ but
$\cW(\s) \not \in \{u,v\}$ for $\s \in [T,t-1]$.
Then $\cB_u \subseteq \cA_u$ and so $\Pr(\cB_u) \le \Pr(\cA_u) $.
It follows from \eqref{ohno}, and $p_uL=O(\pi_S L)$ that
\begin{equation}\label{Buu}
\Pr(\cB_u)  \le p_u(1-(p_u+p_v))^t (1+O(\xi)) .
\end{equation}
However, by contracting $S$ we have that the probability of a first visit to
$\g(S)$ at step $t$ is
\beq{PrES}
\Pr(\cB_u \cup \cB_v)= (1+O(T\pi_S)) (p_u+p_v)(1-(p_u+p_v))^t.
\eeq
From the above and \eqref{Buu}
\begin{eqnarray}
\Pr(\cB_v) & \ge &\Pr(\cB_u \cup \cB_v)-\Pr(\cB_u) \nonumber\\
&\ge& (1-O(\xi))p_v(1-(p_u+p_v))^t.\label{Bvv}
\end{eqnarray}
Using  \eqref{ohno}, \eqref{PrES}, \eqref{Bvv} and $\cE_S=\cB_u \cup \cB_v$  the result follows from
\[
\Pr(\cE_v\mid \cE_S)= \frac{\Pr(\cB_v)}{\Pr(\cB_u \cup \cB_v)}\le
\frac{\Pr(\cA_v)}{\Pr(\cB_u \cup \cB_v)}.
\]
\proofend

\section{Proof of Theorem \ref{th1}(a)}
To apply the lemmas of the previous section we will need to estimate $R_v$
as given by \eqref{RvDef}.

\begin{lemma}\label{lemRv}
If \bP2, \bP3 hold, then in the lazy walk,  for any $v \in V$
\begin{enumerate}[(i)]
\item
\[
 R_v = 2+\frac{2}{d}+O\bfrac{1}{d^{2}}.
\]
\item Suppose $\cW(0)$  is at distance at least 2 from $v$ (resp. at least 3 from $v$). The probability $\cW$ visits $N(v)$ within $L=O(T \log n)$ steps is $P(2,L)=O(1/d)$ (resp. $P(3,L)=O(1/d^{2})$).
\item Let $C\subseteq N(v)$. For a walk starting from $u \in C$, let $R_C$ denote the expected number of returns to $C$ during $T$.
 Then, in the lazy walk,   $R_C=2+O\brac{1/d}$.
\end{enumerate}
\end{lemma}
\proofstart
{\em Proof of (i).}
We write
\[
R_v=1+\sum_{k=1}^{ T}\frac{1}{2^k}+ \sum_{k=0}^{ T-1}\frac{1}{2^k}\sum_{w\in N_G(v)}\frac{1}{2d} \;R(w,T-k-1),
\]
where for $w\in N_G(v)$,  $R(w,\t)$ is the expected number of visits to $v$ in $\t$ steps by $\cW_w$.

For a lower bound,
\[
R(w, \t) \ge \sum_{j=0}^{\t-1} \frac{1}{2^j} \frac{1}{2d}\;R_v=\frac{R_v}{d} \brac{1-\frac{1}{2^{\t}}}.
\]
This is the probability that for some number of  steps  the walk loops at  vertex  $w$, and then moves to $v$,
giving $R_v$ expected returns to $v$. Thus
\[
R_v \ge 2-\frac{1}{2^{T+1}}+ \frac{R_v}{2d}\sum_{k=0}^{ T-1}\frac{1}{2^k}
\brac{1-\frac{1}{2^{T-k-1}}},
\]
so
\[
 R_v \ge 2+\frac{2}{d} +O(1/d^2)-O(T/2^T).
\]
As $T\ge K \log n$ (see \bP1)
we can  assume that $T2^{-T}=O(d^{-2})$.

We next prove we can bound $R(w,T)$ from above by
\beq{EQ0}
R(w,T)\le R_v \brac{\frac{1}{d}+O\bfrac{1}{d^{2}}}.
\eeq
Let $N^i_G(v)$ is the set of vertices at distance $i$ from $v$ in $G$,
and let $R_i^*=\max_{w\in N^i_G(v)}R(w,T)$.
By definition $R(w,T) \le R_1^*$ and
\beq{EQ1}
 R_1^*\leq \sum_{j \ge 0} \brac{\frac12+\frac{\r_2}{2d}}^j\frac{1}{2d}\;R_v+\sum_{j \ge 0} \brac{\frac12+\frac{2\r_2}{2d}}^j\frac{1}{2d}R_1^*+R_3^*.
\eeq
The first summation term counts the case that for some number of  steps  the walk loops at a vertex of $N_G^1(v)$, or moves around in $N_G^1(v)$, which by \bP3 has probability at most $\r_2/2d$. At some point, the walk either moves to $v$,
giving a $R_v$ expected returns, or moves to $N_G^2(v)$.
In the latter case, the second term counts moves back to $N_G^1(v)$, and the third term moves to $N_G^3(v)$,
 giving the $R_3^*$ upper bound.

We next show that $R_3^*=O(1/d^2)$.
Let $\cX$ be
random walk on $\{0,1,\ldots,\r_3\}$, with absorbing barriers at $0,\r_3$,
and transition probabilities for $\cX(i)$  for $0<i<\r_3$  given by
$$\cX(i+1)=\begin{cases}\cX(i)-1&\text{Probability }\frac{\r_2\r_3}{d}\\
\cX(i)&\text{Probability }\frac12\\
\cX(i)+1&\text{Probability }\frac12-\frac{\r_2\r_3}{d}\end{cases}.
$$
Starting $W_z$ from $z \in N^3_G(v)$ and $\cX=\cX_2$ from $j=3$, we can  couple $\cW_z$ and $\cX$ so that $\cX$ is always as close to $0$ as $\cW_z$ is to $v$. Let $u=\cW_z(t)$.
 If $dist(v, u)>\r_3$  then $\cX$ is closer to $v$, where as if
$dist(v, u)\le \r_3$, where $\r_3 \le d^{\e}$, then referring to \bP3, $\nu(v,u) \le \r_2\r_3$.
The probability that $\cW_z(t)$ moves towards $v$ is at most the probability that $\cX$ moves towards $0$.

For a random walk  on $0,1,\ldots,\ell$ starting from $j=0,1,2,\ldots,\ell$
and with probabilities $p,q$ of moving right or left respectively,
it follows from XIV(2.4) of Feller \cite{Feller} that the probability $\p_j$
of the walk visiting 0 before visiting $\ell$ is
\beq{pij}
\p_j=\frac{\xi^j-\xi^{\ell}}{1-\xi^{\ell}}\leq 2\xi^j
\eeq
where $\xi=q/p$. Thus for $\cX$ as given above, $\xi={\r_2\ell}/(d-2\r_2\ell)$.

To finish the proof  of (i),  we choose  $\ell=\r_3=\rdup{d^{\d}}$,
for some $\e/2<\d<\e$.
%and $$\Pi_{T}=\max_{\r_3 \le t \le  T}\Pr(\exists\t,\; t<\t\le T: \cX(\t)=0 \mid %\cX(t)=\r_3).$$
The probability $\pi_3$ that $\cX$ reaches 0 before $\r_3$ is $O(1/d^{3-3\d})=O(1/d^2)$.
Once the walk $\cX$  has reached $\ell=\r_3$,
 we can restart it at $\r_3-1$. The probability it reaches to the origin
 before a return to $\r_3$ is  given by $\pi_{\r_3-1}=O(\xi^{\r_3-1})$.
%As $\Pi_T =O(T \pi_{\r_3-1}) = O( T /d^{\r_3-1})$, and f
From \bP1, $T=O(d^{\r_1} \log n)$, and  we find
\[
R_3^*\le T\p_{\r_3-1}+\p_3= O(\log n \;d^{\r_1+1-\r_3(1-\d)})+O(1/d^{2})=O(1/d^{2}).
\]
For the last inequality, we used $\d>\e/2$ and \bP2 to give
\[
d^{\d} \ge (\log \log n)^{2\d/\e} > \log \log n.
\]
{\em Proof of (ii).}
Let $C=\{v\} \cup N(v)$.
The property \bP3 holds in $G$ for any vertex at distance $\ell \le d^{\e}$
from $v$.
Because moving closer to $C$
implies moving closer to $v$, a vertex
 within distance $\ell$ of $v$ has at most $\r_2 \ell$
neighbours closer to $C$.
If the walk starts at distance 2 from $v$,
it either
loops and/or moves within $N_G^2(v)$, or, conditional on making a transition away from $N_G^2(v)$,
with probability $O(2\r_2/d)$ it moves to $C$, and with probability $1-O(1/d)$
 moves to $N_G^3(v)$.

 Assume the walk starts at  a distance 3 from $v$. We define a graph $\G_{C}$ obtained from
$G$ by contracting the vertices in $C$
to a single vertex $\g_{C}$.
As
explained before Lemma \ref{contract}, we can still use the same  mixing time $T$.
If we replace $v$ by $\g_C$, we can still use the coupling with the random walk $\cX$ on $\{0,1,...,\r_3\}$.
 As moving closer to $\g_C$ means moving closer to $v$, choosing $\r_3=\rdown{d^{\e}}-1$, it follows from \bP3 as outlined above
that the transition probabilities are correct.
By the  argument of part (i), the walk next moves to a distance $\r_3$
from $\g_C$
with probability $1-O(1/d^2)$. After this we use the argument of (i) as before.
In conclusion, for a set $C \subseteq N(v)$ and a walk which moves away from $C$ to a distance 2 from $v$, (resp. distance
3 from $v$)
the probability of a return to $\{v\} \cup N(v)$ within $L$ steps is $O(1/d)$ (resp. $O(1/d^{2})$).

{\em Proof of (iii).}
Let $C \subseteq \{v\} \cup N(v)$. Contract $C$ to $\g_C$ as above.
 We  claim that
$R_{\g_{C}}=2+O\bfrac1{d}$. The 2 comes from the laziness loop at each vertex
and a factor of  $O(\r_2/d)$ comes from possible loops at $\g_{C}$ arising from
$G$-edges inside $C$. If the walk moves to $N_G^2(v)$, then by (ii) the probability of a return to $C$ is $O(1/d)$.
\proofend

\subsection*{Analysis for $t\leq t_{-\e}$}
Recall that $t_{-\e}=(1-\e) n \log d$.
Let $U$ denote the set of vertices unvisited
by the lazy walk in the time interval $[1,2t_{-\e}]$ and
let $U_0$ denote the set of vertices unvisited
by the lazy walk in the time interval $[T,2t_{-\e}]$.
Note that $|U_0\setminus U|\leq T$.
Given Lemma \ref{lem1} below holds,  using
\bP1, \bP2 it follows that $T=o(|U|)$ and thus $|U|=\ooi|U_0|$.

\begin{lemma}\label{lem1} w.h.p.
 \[
 |U_0|\sim \frac{n}{d^{1-\e}} .
 \]
\end{lemma}
\proofstart
Fix a vertex $v$. Corollary \ref{geom} and Remark \ref{remlab}  tell us that
\beq{Probv}
\Pr(v\in U_0)=
\brac{1+O\bfrac{T}{n}}\exp\set{-\frac{2t_{-\e}}{nR_v}+O\bfrac{t_{-\e}}{n^2}} +O(e^{-\Omega(t_{-\e}/T)}).
\eeq
By Lemma \ref{lemRv}, $R_v=2+\frac{2}{d}+O\bfrac1{d^{2}}$.
This gives $\Pr(v\in U_0)\sim d^{1-\e}$ and thus
\[
\E |U_0|\sim \frac{n}{d^{1-\e}}.
\]
Now consider a pair of vertices $v,w$ at distance 5 or more in $G$. Let $\G_{vw}$ be obtained from $G$ by contracting $v,w$ to a single vertex $\g_{vw}$.
%We can couple a random walk $\cW_u'$ on $\G_{vw}$ with a random walk $\cW_u$ on %$G$ up until the time that $\g_{vw}$ is visited. Each choice of the next vertex %has the same probability in both chains until $v$ or $w$ is reached in $G$ and %$\g_{vw}$ is reached in $\G$.
Referring to Lemma \ref{contract} we have
\beq{EQ6}
\Pr(v,w\in U_0)= \Pr(\g_{vw}\in U_0) +O(1/n^3).
\eeq
Working in $\G_{vw}$, it follows more or
less verbatim by using the arguments of Lemma \ref{lemRv}(i) that
$R_{\g_{vw}}=2+\frac2{d}+O\bfrac1{d^{2}}$.
As $v,w$ are sufficiently far apart, only minor
modifications are needed for the analysis of $\cX$. Thus
\beq{2over}
\frac2{R_{\g_{vw}}}=\brac{1+O\bfrac1{d^{2}}}\brac{\frac{1}{R_v}+\frac1{R_w}}.
\eeq
Using $t_{-\e}=(1-\e) n \log d$ in \eqref{Probv} it follows from \eqref{EQ6}
and \eqref{2over} that
\[
\Pr(v,w\in U_0)= \brac{1+ O\bfrac{\log d}{d^{2}}}\;\Pr(v\in U_0)\Pr(w\in U_0) +O(1/n^3).
\]
We prove concentration using the Chebychev inequality. Let $X_{vw}$ be the indicator for $v,w \in U_0$. Let $S$ be  the set of pairs of vertices at distance
at least 5, and let $S'$ be the set of distinct pairs at distance at most 4.
Then
\begin{align*}
\E |U_0|^2 &= \E |U_0| + \sum_{(v,w) \in S} \E X_{vw}+
\sum_{(v,w) \in S'} \E X_{vw}\\
&\le \E |U_0|+ \brac{1+ O\bfrac{\log d}{d^{2}}}\E |U_0|^2+
O(d^4 \E |U_0|).
\end{align*}
It follows from \bP2 that $d^4=o(\E|U_0|)$. Thus for some  $\om$ tending to infinity
\[
\Pr \brac{\mid |U_0|-\E|U_0| \mid \le \frac{\E|U_0|}{\sqrt{\om}}} \le O\bfrac{\om \log d}{d^2} +O\bfrac{
\om d^4}{\E |U_0|}=o(1).
\]
\proofend

%%%%%%%%%central lemma
%%%%%%%%%%%%%%%%%%%%

\begin{lemma}\label{lem2}
 A vertex is {\em bad} if it has fewer than $d^{\e}/2$ neighbours in $U$. Let $B$ denote the set of bad vertices. Then w.h.p.
$|B|\leq ne^{-d^\e/10}$.
\end{lemma}
\proofstart
Fix a vertex $v$ and denote $N_G^1(v)$ by $W=\set{w_1,w_2,\ldots,w_d}$.
Let $X=|W\cap U|$. In the proof of Lemma \ref{lem1} we showed that for a given vertex $x$, $\Pr(x\in U)= \tilde p \sim d^{-(1-\e)}$.
Thus $\E X\sim d^\e$ and if $X$ was distributed as $Bin(d,\tilde p)$ then it would be easy to show that
\beq{1}
\Pr\brac{X\leq \frac12 d^\e}\leq e^{-\Omega(d^\e)}.
\eeq
The bound \eqref{1} is our target. We establish it is true, in spite of $X$ not having a binomial distribution.
For $S\subseteq W$, let $\cA_S=\set{W\cap U=W\setminus S}$, i.e. exactly the vertices  $S$ of $W$ are visited by the walk. So,
\beq{5}
\Pr\brac{X\leq \frac12 d^\e}= \sum_{D=d- d^\e/2}^d \;\sum_{\substack{S\subseteq W\\|S|=D}}\Pr(\cA_S).
\eeq
%Let $\bar{S}=W\setminus S$.
 If $\cA_{S}$ occurs then there is a sequence of times $\bt=(t_0=1\leq t_1<t_2\cdots<t_D\leq t_{D+1}=2t_{-\e})$ and a bijection $f:{S}\to[D]$ such that for $x\in {S}$ there is a first visit to $w_{x}$ at time $t_{f(x)}$.
Let $\cB(S,\bt)$ denote this event.
For a sequence \bt, let $\F(\bt)=\set{i: |t_{i+1}-t_i|\leq L}$, where $L=2KT\log n$. Let $\cT_h=\set{\bt:|\F(\bt)|=h}|$. For $h \ge 0$, let
$$S_h=\sum_{\bt\in \cT_h}\Pr(\cB(S,\bt)).$$
Then,
\beq{6}
\Pr(\cA_S)\leq \sum_{h=0}^D S_h.
\eeq
The main content of the proof of this lemma will be to establish  that
\beq{thething}
\Pr(\cA_S)=O(1)
\brac{e^{-2pt_{-\e}}}^{(d-D)}\brac{1-e^{-2pt_{-\e}}}^D.
\eeq
Given \eqref{5} and \eqref{thething}  we see that
\[
\Pr\brac{X\leq \frac12 d^\e}=O(1) \sum_{D \ge d-d^\e/2}
\binom{d}{D}\brac{e^{-2pt_{-\e}}}^{(d-D)}\brac{1-e^{-2pt_{-\e}}}^D.
\]
The expected value  of $ Bin(d,e^{-2pt_{-\e}})$ is $d^{\e}(1+o(1))$,
so  from the Hoeffding inequality,
\[
\Pr\brac{X\leq \frac12 d^\e}=O\brac{e^{-d^\e/8}}.
\]
Thus, once we prove \eqref{thething}, the lemma  follows from the Markov inequality.

\paragraph{Proof of \ref{thething}.}

We begin with $S_0$.
Our upper bound for $S_0$ will contain some terms that should properly be assigned to some $S_h,h>0$, but this is allowable as we proving an upper bound. We repeat this warning below.
Let
\beq{eqp}
p=\frac{1}{\brac{2+O\bfrac{1}{d}}n},
\eeq
then we have
\begin{multline}\label{2}
S_0\leq D!\sum_{t_1<t_2\cdots<t_D}\brac{\prod_{i=1}^{D} \frac{(1+O(T/n))p}{(1+(d-i+1)p)^{t_i-t_{i-1}}}+ o(e^{-\Omega(\frac{t_i-t_{i-1}}{T}})}\\
\times \brac{\frac{1+O(Td/n)}{(1+(d-D)p)^{2t_{-\e}-t_{D}}}+
o(e^{-\Omega(\frac{t_{-\e}-t_D}{T}})}.
\end{multline}

%%%%%%%%%%%%%%%%%%%%%%%%%%%%%%%%%%%%%%%%%%%%%%%%%%
%%%%%%%%%%%%%%%%%%%%%%%%%%%%%%%%%%%%%%%%%%%%%%%%%%%

{\bf Proof of \eqref{2}:} Assume for the moment that $S=\set{w_1,\ldots,w_D}$ and that $f(w_i)=i$ for $i=1,2,\ldots,D$. Let  $A_i=\set{w_i,w_{i+1},\ldots,w_D}$ for $i=1,2,\ldots,D$. We assign times $t_1,t_2,\ldots,t_D$ to ${S}$ in $D!$ ways. Now consider a term
\beq{PSI}
\Psi_i=\frac{(1+O(T/n))p}{(1+(d-i+1)p)^{t_i-t_{i-1}}}
+o(e^{-\Omega((t_i-t_{i-1})/T)}).
\eeq
We claim this is an estimate of the probability
there are no visits to $w_{i},\ldots,w_D$ during $[t_{i-1}+T,t_i-1]$ followed by
 a first visit to $w_i$ at $t_i$.
If so, it is also  an upper bound for the probability  there is
 no visit to $w_{i},\ldots,w_D$ during $[t_{i-1}+1,t_i-1]$ followed by a visit to $w_i$ at $t_i$.
 This bound hold regardless of the first $t_{i-1}$ steps of the walk. In fact this bound allows for visits to $w_i,w_{i+1},\ldots,w_D$ during the time interval $[t_{i-1}+1,t_{i-1}+T-1]$, but this is allowable as $\Psi_i$ is an upper bound. Thus some terms properly attributed to $S_h,h>0$ are overcounted.

To prove \eqref{PSI},  define a graph $\G_{A_i}$ obtained from $G$ by contracting the vertices in $A_i$
to a single vertex $\g_{A_i}$. The mixing time  $T$ does not increase, as explained above Lemma \ref{contract}. We also have $R_{\g_{A_i}}\leq2+O\bfrac1{d}$. For this, we again follow the proof of Lemma \ref{lemRv}. The 2 comes from the laziness loop at each vertex and the $O\bfrac{1}{d}$ comes from possible loops at $\g_{A_i}$ arising from cases where there are $G$-edges inside $A_i$. We apply the same argument as in Lemma \ref{lemRv}. We can use the random walk $\cX$ because a vertex $z\neq \g_{A_i}$ and within distance $\r_3-1$ of $\g_{A_i}$ has at most $\r_2\r_3$ neighbours closer to $\g_{A_i}$. This is because moving closer to $\g_{A_i}$ implies moving closer to $A_i  \subseteq W$, and hence to $v$. Apply \bP3 to $\{z,v\}$.

By Lemma \ref{First},  the probability $t_{i}$ is the time of a first visit to $\g_{A_i}$ in $[t_{i-1}+T,t_i]$ can be expressed as $(d-i+1)\Psi_i$.
Given a first visit has been made to $A_i$, we need  the
probability that this first visit was made to a given $v \in A_i$. Lemma \ref{whichvx} gives the answer. The $p_{w_j}, j=i,...,D$ used in Lemma \ref{whichvx} are given by \eqref{eqp}. This establishes \eqref{PSI}.

The final term in \eqref{2}, given by $\frac{1+O(Td/n)}{(1+(d-D)p)^{2t_{-\e}-t_{D}}}+o(e^{-\Omega((t_{-\e}-t_D)/T)})$ bounds the probability that the vertices in $\set{w_{D+1},\ldots,w_d}$ are not visited in the interval $[t_D,2t_{-\e}]$. We use the first part of the argument for $\Psi_i$ to validate this.\\
{\bf End of proof of \eqref{2}.}

\medskip
The next step is to evaluate \eqref{2}. Considering \eqref{PSI},
the term $\frac{p}{(1+(d-i+1)p)^{t_i-t_{i-1}}}=\Omega((1/n)e^{(t_i-t_{i-1})/n})$, whereas  the term
$o(e^{-\Omega((t_i-t_{i-1})/T})=o(e^{(t_i-t_{i-1})/T})$. As $t_i-t_{i-1} \ge L=KT \log n$ the
latter term can be absorbed into the $O(d^{-1})$ in the definition of
$p$. Furthermore,
$$\frac{1}{1+(d-i+1)p}=\exp\set{-(d-i+1)p+O\bfrac{d^2}{n^2}}.$$
Noting that
\[
\sum_{i=1}^{D+1}(d-i+1)(t_i-t_{i-1})=(d-D)t_{D+1}+(t_1+\cdots+t_D),
\]
we can write
\begin{align}
S_0&\leq
2D!p^D\sum_{t_1<t_2\cdots<t_D}\exp\set{-p\sum_{i=1}^{D+1}(d-i+1)(t_i-t_{i-1})}\nonumber\\
&=2D!p^De^{-2(d-D)pt_{-\e}}\sum_{t_1<t_2\cdots<t_D}\exp\set{-p\sum_{i=1}^{D}t_i}\nonumber\\
&\leq 2e^{-2(d-D)pt_{-\e}}\brac{p\sum_{t=1}^{2t_{-\e}}e^{-pt}}^D\nonumber\\
&\leq3e^{-2(d-D)pt_{-\e}}\brac{p\int_{t=0}^{2t_{-\e}}e^{-pt}dt}^D\nonumber\\
&=3e^{-2(d-D)pt_{-\e}}(1-e^{-2pt_{-\e}})^D\label{binoval}
\end{align}
We next show  that $S_1,S_2.\ldots,S_D$ are not much larger in total than
$S_0$.

We say a  visit to vertex $u$ is {\em $T$-distinct}, if it occurs at least $T$ steps after a previous $T$-distinct visit, or from the start of the walk. Thus if $\cW(t)=u$, and this visit is $T$-distinct, the next $T$-distinct visit to $u$ will be at the first step $s \ge t+T$ such that $\cW(s)=u$.
Once a $T$-distinct visit has  taken place, several {\em secondary visits} to the vertex $u$ may occur within the next $T-1$ steps, and thus before the next $T$-distinct visit. We will consider such secondary visits separately in our proof.

We consider the case  $t_i-t_{i-1} \le L$ in two parts, namely $t_i-t_{i-1} < T$, and $T \le t_i-t_{i-1} \le L$.
The first case is for {\em secondary visits},  and the second case {\em close (together) visits}.
These require a separate analysis.

Given $\bt=(t_1,\ldots,t_D)$ for arbitrary $D \le d$, let $Z\ge D-k$ be an upper bound on the total number of secondary visits to $W=N(v)$
occurring as a result of $k\le D$ first visits to $W$ which are $T$-distinct.
Let $N_2(v)$ denote the set of vertices at distance 2 from $v$. Then
\[
Z(\bt)= N_1+\cdots+N_k
\]
where $N_i$ are the number of secondary visits to $W=N(v)$ (i.e. returns to $W$ via $\{v\} \cup N_2(v)$) which occur during $[t_i,t_i+T]$, $i=1,...,k$.

The values of $N_i$ are independent and   geometrically distributed with failure probability $O(1/d)$.
From $W=N(v)$ the particle  moves to $\{v\} \cup N(v)$  with probability $O(1/d)$, (this  follows from \bP3).
Otherwise the particle moves to distance 2 away from $v$ with probability $1-O(1/d)$, and we can use the value of $P(2,T)=O(1/d)$ from
 Lemma \ref{lemRv}(ii).
For any $D \le d$, the probability $\widehat P(\ell)$ of at least $\ell$ secondary visits is
\[
\widehat P(\ell) = \binom{D+\ell-1}{\ell}\bfrac{O(1)}{d}^\ell \le \bfrac{O(1)D}{\ell d}^\ell\le \bfrac{O(1)}{\ell}^\ell=e^{-\Th(\e d^{\e} \log d )},
\]
 on choosing $\ell =d^{\e}/100$. Provided $\e \gg 1/\log d$, the probability of at least $d^{\e}/100$ secondary visits to $W$ is $o(e^{-d^\e})$.

We next consider close together visits. For convenience, replace $D$ by $D'=D-Z$ i.e. remove any entries in $\bt$ corresponding to secondary
visits.
Let $h$ count those $T$-distinct first visits which are close together i.e. $T \le t_i-t_{i-1} \le  L$.
After $t \ge T$ steps, the distribution of the walk is close to stationary, so  the probability that the walk is within distance 2 of vertex
$v$ is $O(d^2/n)$. If the walk is at least distance 3 from $v$,  by
 Lemma \ref{lemRv}(ii)  the probability of a visit to $W=N(v)$ in $L$ steps is at most $P(3,L)=O(1/d^{2})$.
It follows that, independently of any previous ones, each close visit has probability $O(d^2/n)+O(1/d^{2})=O(1/d^{2})$,
assuming $d=o(n^{1/4})$ (see \bP2).

To bound $S_h$ we note that the remaining $k=D-h$ first visits are \lq well spaced' i.e. $L \le t_i-t_{i-1}$.
There are $\binom{D-1}{h}$ ways to assign the $h$ \lq close together' events to the $k=D-h$ \lq well spaced' ones.
To do so, we choose an allocation $n_1,n_2,\ldots,n_k\geq 0$ such that $n_1+n_2+\cdots+n_k=h$.

Note that $S_0=S_0(D)$ so changing $D$ to $D-h$, for $h\ge 1$, from \eqref{binoval} we have
\beq{sht}
S_h(D) \le  S_0(D-h) \binom{D-1}{h} \bfrac{O(1)}{d^2}^h \le S_0(D) \bfrac{e^{2pt_{-\e}}}{1-e^{-2pt_{-\e}}}^h \bfrac{O(D)}{hd^2}^h\le
S_0(D) \bfrac{O( d^{-\e})}{h}^h.
\eeq
The value of $p$ is from \eqref{eqp}, and $t_{-\e}=(1-\e)n \log d$.
Inequality \eqref{sht}, along with \eqref{binoval} completes the proof of \eqref{thething},
and the lemma follows.
\proofend

We can now easily show that w.h.p. at time $2t_{-\e}$, there is a
component of size much larger than $\log n$.
\begin{lemma}\label{lem3}
W.h.p. the graph induced by unvisited vertices contains a component of
size at least $e^{\Omega(d^{\e/2})}$.
\end{lemma}
\proofstart
Let $n_0=\frac{n}{5(ed^{1-\e/2})^{d^{\e/2}}d^{1-\e}}$. We begin by greedily choosing $v_1,v_2,\ldots,v_{n_0}\in U$ such that
$v_i,v_j$ are at distance greater than $d^{\e/2}$. This is easily done,
because there are
$1+\binom{d}{1}+\binom{d}{2}+\cdots+\binom{d}{d^{\e/2}}<2\binom{d}{d^{\e/2}}\leq
2(ed^{1-\e/2})^{d^{\e/2}}$
vertices within distance $d^{\e/2}$ of any given vertex. Having chosen
$v_1,v_2,\ldots,v_k,\,k\leq n_0$, there will w.h.p. be at least
$\frac{n}{2d^{1-\e}}-2k(ed^{1-\e/2})^{d^{\e/2}}>0$ choices for $v_{k+1}$. For each $i$ let
$V_i$ denote the set of vertices within distance $d^{\e/2}$ of $v_i$. The
$V_i$ are disjoint and so from Lemma \ref{lem2} there are w.h.p. at least $n_0-ne^{-d^\e/10}>0$ indices $i$
such that $V_i\cap B=\emptyset$.

Choose $i$ such that $V_i\cap B=\emptyset$. From $v_i$ we can do
breadth first search, but only including vertices in $U$. If $L_r$
denotes the $r$th level of this search where $L_0=\set{v_i}$ then we
see that $|L_{r+1}|\geq \frac{d^\e|L_r|}{2\r_2(r+1)}$. Thus $V_i$ contains
a component of size at least
$$\sum_{i=0}^{d^{\e/2}/2}\binom{d^\e/2}{i}\frac{1}{(2\r_2)^i}=e^{\Omega(d^{\e/2})}.$$
\proofend
\section{Proof of Theorem \ref{th1}(b)}
Let
$$s=\frac{2\log n}{\e \log d}=o(\log n).$$
We will show that w.h.p. there is no component of size $s$ or more at time $t\geq 2t_{+\e}$ in $\G(t)$, with respect to the lazy walk.

\begin{lemma}\label{cl}
For $v\in V$ there are at most $(ed)^{s-1}$ sets $S$ such that (i) $v\in S$, (ii) $|S|=s$ and (iii) $G[S]$ is connected.
\end{lemma}
\proofstart The number of such sets is bounded by the number of distinct $s$-vertex
trees which are rooted at $v$. This in turn is bounded by the number of distinct $d$-ary rooted trees
with $s$ vertices. This is equal to $\binom{ds}{s}/ ((d-1)s+1)$, see Knuth \cite{Knuth}.
\proofend

We fix a set $S$ of size $s$ that induces a connected subgraph of
$G$. To estimate the probability that $S$ is unvisited at time $t\geq
2t_{+\e}$ we contract $S$ to a vertex $\g_S$ as in the proofs of
Lemmas \ref{lem1} and \ref{lem2}. We need to estimate the probability that $\g_S$ is unvisited by a lazy random walk on the associated graph $\G_S$ during the time interval $[T,2t_{+\e}]$. For this we need to prove
\begin{lemma}\label{lemRg}
$R_{\g_S}=2+o(1)$.
\end{lemma}
\proofstart
Let $e(S)$ denote the number of edges contained in $S$. It follows from \bP4 that $e(S)=o(ds)$. This means that $\g_S$ has degree $ds$, of  which $o(ds)$ comes from loops associated with internal edges of $S$. It then follows that when the walk on $\G_S$ is at $\g_S$ then it leaves $\g_S$ with probability $\frac12-o(1)$. It is then straightforward to use the argument of Lemma \ref{lemRv} to finish the proof of the lemma.
\proofend

Using Lemma \ref{cl} and Lemma \ref{lemRg} we see that if
$p_\g=\frac{(1+o(1))s}{2}$
then
\begin{align*}
\Pr( \text{there exists a component of size }s)&\leq n(ed)^{s-1}\brac{\frac{1+O(Ts/n))}{(1+p_\g)^{2t_{+\e}}}+O(T^2se^{-\Omega(t_{+\e}/T)})}\\
&\leq 2n(ed\cdot e^{-(1-o(1))(1+\e)\log d})^s\\
&\leq 2nd^{-2\e s/3}=o(1).
\end{align*}
\proofend

%%%%%%%%%%%%%%%%%%%%%%%%%%%%%%%%%%%%%%%%%%%%%%%%%%%%%%%%%%%%%%%%%%%%

\end{document}